\documentclass[amstex,12pt, amssymb]{article}

\usepackage{mathtext}
\usepackage[cp1251]{inputenc}
\usepackage[T2A]{fontenc}
\usepackage[dvips]{graphicx}
\usepackage{amsmath}
\usepackage{amssymb}
\usepackage{amsxtra}
\usepackage{latexsym}
\usepackage{ifthen}

\newcounter{lemma}[section]

\newcounter{corollary}[section]

\newcounter{remark}[section]

\newcounter{theorem}[section]

\newcounter{proposition}[section]

\newcounter{example}

\numberwithin{equation}{section}

\begin{document}

\markboth{E~.SEVOST'YANOV, V.~TARGONSKII}{\centerline{ON CONVERGENCE
OF HOMEOMORPHISMS...}}

\def\cc{\setcounter{equation}{0}
\setcounter{figure}{0}\setcounter{table}{0}}

\overfullrule=0pt

%\normalsize\large

\author{EVGENY SEVOST'YANOV, VALERY TARGONSKII}

\title{
{\bf ON CONVERGENCE OF HOMEOMORPHISMS WITH INVERSE MODULUS
INEQUALITY}}

\date{\today}
\maketitle

%\large
\begin{abstract}
We have studied homeomorphisms that satisfy the Poletsky-type
inverse inequality in the domain of the Euclidean space. It is
proved that the uniform limit of the family of such homeomorphisms
is either a homeomorphism into the Euclidean space, or a constant in
the extended Euclidean space.
\end{abstract}



\section{Introduction}

This paper is devoted to the study of mappings with bounded and
finite distortion, see, e.g., \cite{Cr$_1$}--\cite{Cr$_2$},
\cite{MRV$_1$}, \cite{MRSY}, \cite{Vu} and~\cite{Va}. It is well
known that the locally uniform limit of quasiconformal mappings is a
homeomorphism, or a constant, see, e.g., \cite[Theorems~21.9,
21.11]{Va}. This fact is true not only for quasiconformal mappings,
but also in broader classes of mappings that satisfy modulus
conditions. In particular, the first co-author together with
V.~Ryazanov proved that the specified property holds for the
so-called ring $Q$-homeomorphisms under certain conditions regarding
the function $Q,$ see, e.g., \cite[Theorems~4.1 and 4.2]{RS}. This
result was generalized by M.~Cristea for more general classes of
mappings and somewhat more general conditions on $Q,$ and for the
so-called weighted modulus conditions
(see~\cite[Theorem~1]{Cr$_2$}). In this manuscript, we will show the
validity of a similar statement for maps with an inverse modulus
condition, i.e., maps inverse to ring $Q$-homeomorphisms. It should
be noted that this statement does not follow from the previously
obtained results, because the image domain of under homeomorphisms
may be variable. In particular, if we make the transition from
direct mappings to inverses, we will not get a class of mappings,
defined in a single domain. In the corresponding previous results,
the assumption that all mappings are defined in a single domain is
essential.

\medskip Below $dm(x)$ denotes the element of the Lebesgue
measure in ${\Bbb R}^n.$ Everywhere further the boundary $\partial A
$ of the set $A$ and the closure $\overline{A}$ should be understood
in the sense of the extended Euclidean space $\overline{{\Bbb
R}^n}.$ Recall that, a Borel function $\rho:{\Bbb R}^n\,\rightarrow
[0,\infty] $ is called {\it admissible} for the family $\Gamma$ of
paths $\gamma$ in ${\Bbb R}^n,$ if the relation
\begin{equation}\label{eq1.4}
\int\limits_{\gamma}\rho (x)\, |dx|\geqslant 1
\end{equation}
holds for all (locally rectifiable) paths $ \gamma \in \Gamma.$ In
this case, we write: $\rho \in {\rm adm} \,\Gamma .$  The {\it
modulus} of $\Gamma $ is defined by the equality
\begin{equation}\label{eq1.3gl0}
M(\Gamma)=\inf\limits_{\rho \in \,{\rm adm}\,\Gamma}
\int\limits_{{\Bbb R}^n} \rho^n (x)\,dm(x)\,.
\end{equation}
Let $y_0\in {\Bbb R}^n,$ $0<r_1<r_2<\infty$ and
\begin{equation}\label{eq1**}
A=A(y_0, r_1,r_2)=\left\{ y\,\in\,{\Bbb R}^n:
r_1<|y-y_0|<r_2\right\}\,.\end{equation}
Given $x_0\in{\Bbb R}^n,$ we put
$$B(x_0, r)=\{x\in {\Bbb R}^n: |x-x_0|<r\}\,, \quad {\Bbb B}^n=B(0, 1)\,,$$
$$S(x_0,r) = \{
x\,\in\,{\Bbb R}^n : |x-x_0|=r\}\,. $$
A mapping $f: D \rightarrow{\Bbb R}^n$ is called {\it discrete} if
the pre-image $\{f^{-1}\left(y\right)\}$ of any point $y\,\in\,{\Bbb
R}^n$ consists of isolated points, and {\it open} if the image of
any open set $U\subset D$ is an open set in ${\Bbb R}^n.$

Given sets $E,$ $F\subset\overline{{\Bbb R}^n}$ and a domain
$D\subset {\Bbb R}^n$ we denote by $\Gamma(E,F,D)$ the family of all
paths $\gamma:[a,b]\rightarrow \overline{{\Bbb R}^n}$ such that
$\gamma(a)\in E,\gamma(b)\in\,F$ and $\gamma(t)\in D$ for $t \in (a,
b).$ Given a mapping $f:D\rightarrow \overline{{\Bbb R}^n},$ a point
$y_0\in {\Bbb R}^n,$ and $0<r_1<r_2<r_0=\sup\limits_{y\in
f(D)}|y-y_0|,$ we denote by $\Gamma_f(y_0, r_1, r_2)$ a family of
all paths $\gamma$ in $D$ such that $f(\gamma)\in \Gamma(S(y_0,
r_1), S(y_0, r_2), A(y_0,r_1,r_2)).$ Let $Q:{\Bbb R}^n\rightarrow
[0, \infty]$ be a Lebesgue measurable function. We say that {\it $f$
satisfies the inverse Poletsky inequality at a point $y_0\in {\Bbb
R}^n$} if the relation
\begin{equation}\label{eq2*A}
M(\Gamma_f(y_0, r_1, r_2))\leqslant \int\limits_{A(y_0,r_1,r_2)\cap
f(D)} Q(y)\cdot \eta^{n}(|y-y_0|)\, dm(y)
\end{equation}
holds for any Lebesgue measurable function $\eta:
(r_1,r_2)\rightarrow [0,\infty ]$ such that
\begin{equation}\label{eqA2}
\int\limits_{r_1}^{r_2}\eta(r)\, dr\geqslant 1\,.
\end{equation}
The definition of the relation~(\ref{eq2*A}) at the point
$y_0=\infty$ may be given by the using of the inversion
$\psi(y)=\frac{y}{|y|^2}$ at the origin.

\medskip
Note that conformal mappings preserve the modulus of families of
paths, so that we may write
$$M(\Gamma)=M(f(\Gamma))\,.$$

\medskip
Set
\begin{equation}\label{eq12}
q_{y_0}(r)=\frac{1}{\omega_{n-1}r^{n-1}}\int\limits_{S(y_0,
r)}Q(y)\,d\mathcal{H}^{n-1}(y)\,, \end{equation}
and $\omega_{n-1}$ denotes the area of the unit sphere ${\Bbb
S}^{n-1}$ in ${\Bbb R}^n.$

\medskip
We say that a function ${\varphi}:D\rightarrow{\Bbb R}$ has a {\it
finite mean oscillation} at a point $x_0\in D,$ write $\varphi\in
FMO(x_0),$ if
$$\limsup\limits_{\varepsilon\rightarrow
0}\frac{1}{\Omega_n\varepsilon^n}\int\limits_{B( x_0,\,\varepsilon)}
|{\varphi}(x)-\overline{{\varphi}}_{\varepsilon}|\ dm(x)<\infty\,,
$$
where $\overline{{\varphi}}_{\varepsilon}=\frac{1}
{\Omega_n\varepsilon^n}\int\limits_{B(x_0,\,\varepsilon)}
{\varphi}(x) \,dm(x)$ and $\Omega_n$ is the volume of the unit ball
${\Bbb B}^n$ in ${\Bbb R}^n.$
We also say that a function ${\varphi}:D\rightarrow{\Bbb R}$ has a
finite mean oscillation at $A\subset \overline{D},$ write
${\varphi}\in FMO(A),$ if ${\varphi}$ has a finite mean oscillation
at any point $x_0\in A.$ Let $h$ be a chordal metric in
$\overline{{\Bbb R}^n},$
$$h(x,\infty)=\frac{1}{\sqrt{1+{|x|}^2}}\,,$$
\begin{equation}\label{eq3C}
h(x,y)=\frac{|x-y|}{\sqrt{1+{|x|}^2} \sqrt{1+{|y|}^2}}\qquad x\ne
\infty\ne y\,.
\end{equation}
and let $h(E):=\sup\limits_{x,y\in E}\,h(x,y)$ be a chordal diameter
of a set~$E\subset \overline{{\Bbb R}^n}$ (see, e.g.,
\cite[Definition~12.1]{Va}).

\medskip
\begin{theorem}\label{th1}
{\it Let $D$ be a domain in ${\Bbb R}^n,$ $n\geqslant 2,$ and let
$f_m:D\rightarrow {\Bbb R}^n,$ $m=1,2,\ldots ,$ be a sequence of
homeomorphisms that converges to some mapping $f:D\rightarrow
\overline{{\Bbb R}^n}$ locally uniformly in $D$ by the metric $h,$
and satisfy the relations~(\ref{eq2*A})--(\ref{eqA2}) in each point
$y_0\in \overline{{\Bbb R}^n}.$ Assume that, one of two conditions
holds:

\medskip
1) $Q\in FMO(\overline{{\Bbb R}^n});$

\medskip
2) for any $y_0\in \overline{{\Bbb R}^n}$ there exists
$\delta(y_0)>0$ such that
\begin{equation}\label{eq5B}
\int\limits_{\varepsilon}^{\delta(y_0)}
\frac{dt}{tq_{y_0}^{\frac{1}{n-1}}(t)}<\infty, \qquad
\int\limits_{0}^{\delta(y_0)}
\frac{dt}{tq_{y_0}^{\frac{1}{n-1}}(t)}=\infty
\end{equation}
for sufficiently small $\varepsilon>0.$ Then $f$ is either a
homeomorphism $f:D\rightarrow {\Bbb R}^n,$ or a constant
$c\in\overline{{\Bbb R}^n}.$ }
\end{theorem}

\medskip
Here the conditions mentioned above for $y_0=\infty$ must be
understood as conditions for the function
$\widetilde{Q}(y):=Q(y/|y|^2)$ at the origin.

\section{Preliminaries}

Let $D\subset {\Bbb R}^n,$ $f:D\rightarrow {\Bbb R}^n$ be a discrete
open mapping, $\beta: [a,\,b)\rightarrow {\Bbb R}^n$ be a path, and
$x\in\,f^{\,-1}(\beta(a)).$ A path $\alpha: [a,\,c)\rightarrow D$ is
called a {\it maximal $f$-lifting} of $\beta$ starting at $x,$ if
$(1)\quad \alpha(a)=x\,;$ $(2)\quad f\circ\alpha=\beta|_{[a,\,c)};$
$(3)$\quad for $c<c^{\prime}\leqslant b,$ there is no a path
$\alpha^{\prime}: [a,\,c^{\prime})\rightarrow D$ such that
$\alpha=\alpha^{\prime}|_{[a,\,c)}$ and $f\circ
\alpha^{\,\prime}=\beta|_{[a,\,c^{\prime})}.$ If $\beta:[a,
b)\rightarrow\overline{{\Bbb R}^n}$ is a path and if
$C\subset\overline{{\Bbb R}^n},$ we say that $\beta\rightarrow C$ as
$t\rightarrow b,$ if the spherical distance $h(\beta(t),
C)\rightarrow 0$ as $t\rightarrow b$ (see
\cite[section~3.11]{MRV$_2$}), where $h(\beta(t),
C)=\inf\limits_{x\in C}h(\beta(t), x).$ The following assertion
holds (see~\cite[Lemma~3.12]{MRV$_2$}).

\medskip
\begin{proposition}\label{pr3}
{\it Let $f:D\rightarrow {\Bbb R}^n,$ $n\geqslant 2,$ be an open
discrete mapping, let $x_0\in D,$ and let $\beta: [a,\,b)\rightarrow
{\Bbb R}^n$ be a path such that $\beta(a)=f(x_0)$ and such that
either $\lim\limits_{t\rightarrow b}\beta(t)$ exists, or
$\beta(t)\rightarrow \partial f(D)$ as $t\rightarrow b.$ Then
$\beta$ has a maximal $f$-lifting $\alpha: [a,\,c)\rightarrow D$
starting at $x_0.$ If $\alpha(t)\rightarrow x_1\in D$ as
$t\rightarrow c,$ then $c=b$ and $f(x_1)=\lim\limits_{t\rightarrow
b}\beta(t).$ Otherwise $\alpha(t)\rightarrow \partial D$ as
$t\rightarrow c.$}
\end{proposition}

For a domain $D\subset {\Bbb R}^n,$ $n\geqslant 2,$ and a Lebesgue
measurable function $Q:{\Bbb R}^n\rightarrow [0, \infty],$
$Q(y)\equiv 0$ for $y\in{\Bbb R}^n\setminus f(D),$ we denote by
$\frak{F}_Q(D)$ the family of all open discrete mappings
$f:D\rightarrow {\Bbb R}^n$ such that
relations~(\ref{eq2*A})--(\ref{eqA2}) hold for each point $y_0\in
f(D).$ The following result holds (see~\cite[Theorem~1.1]{SSD}).

\medskip
A domain $R$ in $\overline{{\Bbb R}^n},$ $n\geqslant 2,$ is called
{\it a ring,} if $\overline{{\Bbb R}^n}\setminus R$ consists of
exactly two components $E$ and $F.$ In this case, we write: $R=R(E,
F).$ The following statement is true, see \cite[ratio~(7.29)]{MRSY}.

\begin{proposition}\label{pr2}{\it\, If $R=R(E, F)$ is a ring, then
$$M(\Gamma(E, F, \overline{{\Bbb R}^n}))\geqslant\frac{\omega_{n-1}}{\left(
\log\frac{2\lambda^2_n}{h(E)h(F)}\right)^{n-1}}\,,$$
where $\lambda_n \in[4,2e^{n-1}),$ $\lambda_2=4$ and
$\lambda_n^{1/n} \rightarrow e$ as $n\rightarrow \infty,$ and $h(E)$
denotes the chordal diameter of the set $E,$
$h(E):=\sup\limits_{x,y\in E}\,h(x,y).$}
\end{proposition}

\medskip
In accordance with~\cite{GM}, a domain $D$ in ${\Bbb R}^n$ is called
a {\it quasiextremal distance domain} (a $QED$-{\it domain for
short)} if
\begin{equation}\label{eq4***}
M(\Gamma(E, F, {\Bbb R}^n))\leqslant  A\cdot M(\Gamma(E, F, D))
\end{equation}
for some finite number $A\geqslant 1$ and all continua $E$ and $F$
in $D$.

\medskip
Recall the following statement, see~\cite[Theorem~3.1]{RS}.

\medskip
\begin{proposition}\label{pr4}
{\it Let $D$ be a domain in $\overline{{\Bbb R}^n},$ $n\geqslant 2,$
and let $f_m,$ $m=1, 2, \ldots ,$ be a sequence of homeomorphisms of
$D$ into $\overline{{\Bbb R}^n}$ converging locally uniformly to a
discrete mapping $f:D\rightarrow \overline{{\Bbb R}^n}$ with respect
to the spherical (chordal) metric. Then $f$ is a homeomorphism of
$D$ into $\overline{{\Bbb R}^n}.$}
\end{proposition}

\section{Main Lemmas}

\begin{lemma}\label{lem1}
{\it Let $D$ be a $QED$-domain in ${\Bbb R}^n,$ $n\geqslant 2,$ and
let $f:D\rightarrow {\Bbb R}^n$ be a homeomorphism satisfying the
relations~(\ref{eq2*A})--(\ref{eqA2}) at some $y_0\in {\Bbb R}^n.$
Let $\varepsilon_1>0$ be such that $B(x_0, \varepsilon_1)\subset D,$
let $x\in B(x_0, \varepsilon_1),$ let $\overline{B(z_0,
\varepsilon_2)}$ be a closed ball in $D,$ and let $\varepsilon_0>0$
be such that
\begin{equation}\label{eq1E}
\varepsilon:=|f(x)-f(x_0)|+|f(x_0)-y_0|<\varepsilon_0\,,\qquad
f(\overline{B(z_0, \varepsilon_2)})\cap\overline{B(y_0,
\varepsilon_0)}=\varnothing\,.
\end{equation}
Assume that, there is a Lebesgue measurable function $\psi:(0,
\varepsilon_0)\rightarrow (0,\infty)$ and a constant $c_3>0$ such
that
\begin{equation}\label{eq7***} 0<I(\varepsilon,
\varepsilon_0):=\int\limits_{\varepsilon}^{\varepsilon_0}\psi(t)\,dt
< \infty\,,
\end{equation}
while there exits a function $\alpha=\alpha(\varepsilon,
\varepsilon_0)>0$ such that
\begin{equation} \label{eq3.7.2}
\int\limits_{A(y_0, \varepsilon, \varepsilon_0)}
Q(y)\cdot\psi^{\,n}(|y-y_0|)\,dm(x)= \alpha(\varepsilon,
\varepsilon_0)\cdot I^n(\varepsilon, \varepsilon_0)\,,\end{equation}
where $A(y_0, \varepsilon, \varepsilon_0)$ is defined in
(\ref{eq1**}). Then
\begin{equation}\label{eq1}
|x-x_0|\leqslant \frac{2\lambda^2_n}{c_1\cdot h(\overline{B(z_0,
\varepsilon_2)})}\cdot\exp\left\{-\frac{c_2\omega_{n-1}}{\alpha(|f(x)-f(x_0)|+|y_0-f(x_0)|,
\varepsilon_0)}\right\}\,,
\end{equation}
where $c_1:=\frac{1}{1+\varepsilon^2_1}$ and $c_2$ is a constant
from the definition of $QED$-domain for $D,$ i.e., $c_2:=A$
in~(\ref{eq4***}). }
\end{lemma}

\medskip
{\it Proof of Lemma~\ref{lem1}.} Let $x\in B(x_0, \varepsilon_1).$
Let us join the points $f(x)$ and $f(x_0)$ by the segment $I,$
$I=I(t)=f(x_0)+(f(x)-f(x_0))t,$ $t\in [0, 1].$ Let $\alpha:[0,
c)\rightarrow D$ be a maximal $f$-lifting of $I$ starting at $x_0.$
By Proposition~\ref{pr3} this lifting is well-defined and either one
of the following situations holds: $\alpha(t)\rightarrow x_1\in
B(x_0, \varepsilon_1)$ as $t\rightarrow c-0$ (in this case, $c=1$
and $f(x_1)=f(x)$), or $\alpha(t)\rightarrow S(x_0, \varepsilon_1)$
as $t\rightarrow c.$ In the first situation, $x=x_1$ because $f$ is
a homeomorphism. Choose $\Delta>0$ such that
$|x-x_0|<\Delta<\varepsilon_1,$ and let $t_0=\sup\limits_{t\in [0,
1]: \alpha(t)\in B(0, \Delta)}t.$ Observe that,
$\alpha_1:=\alpha|_{[0, t_0]}$ is a closed Jordan path because
$\alpha_1(t)=f^{\,-1}(I(t))$ and $f$ is a homeomorphism. In
particular, $|\alpha_1|$ is a continuum. If $|\alpha_1|\subset
B(x_0, \Delta),$ then $t_0=1$ and $\alpha_1=\alpha,$
$\alpha_1(1)=\alpha(1)=x.$ Otherwise, $\alpha(t_0)\in S(x_0,
\Delta).$ Thus,
\begin{equation}\label{eq2A}
{\rm diam\,}|\alpha_1|\geqslant \min\{\Delta\,, |x-x_0|\}>|x-x_0|\,.
\end{equation}
Let $\overline{B(z_0, \varepsilon_2)}\cap |\alpha_1|=\varnothing.$
Let us to prove that $R=R(\overline{B(z_0, \varepsilon_2)},
|\alpha_1|)$ is a ring domain.

\medskip
Indeed, since $\alpha_1$ is a Jordan path, it does not split ${\Bbb
R}^n$ for $n\geqslant 3,$ because $|\alpha_1|$ has a topological
dimension~1 (see~\cite[Theorem~III 2.3]{HW} and
~\cite[Corollary~1.5.IV]{HW}). Now, any points $x_1, x_2\in {\Bbb
R}^n\setminus (\overline{B(z_0, \varepsilon_2)}\cup |\alpha_1|)$ may
be joined by a path $\gamma:[0, 1]\rightarrow {\Bbb R}^n,$
$\gamma(0)=x_1,$ $\gamma(1)=x_2,$ in ${\Bbb R}^n\setminus
|\alpha_1|.$

\medskip
Let us to show that the same is true for $n=2.$ Join the points
$x_1, x_2\in {\Bbb R}^n\setminus (\overline{B(z_0,
\varepsilon_2)}\cup |\alpha_1|)$ by some path
$\widetilde{\gamma}:[0, 1]\rightarrow {\Bbb R}^n,$
$\widetilde{\gamma}(0)=x_1,$ $\widetilde{\gamma}(1)=x_2,$ in ${\Bbb
R}^n.$ If $\widetilde{\gamma}\cap|\alpha_1|=\varnothing,$ it is
nothing to prove. Otherwise, due to Antoine's theorem on the absence
of wild arcs (see~\cite[Theorem~II.4.3]{Keld}), there exists a
homeomorphism $\varphi:{\Bbb R}^2\rightarrow {\Bbb R}^2,$ which maps
$\alpha_1$ onto some segment $I.$ Let $\Pi$ be an open rectangular
two of edges of which are parallel to $I,$ and two of which are
perpendicular to $I,$ while $I\subset \Pi.$ Reducing $\Pi,$ we also
may assume that $\widetilde{\varphi}(x_1)\not\in \Pi$ and
$\widetilde{\varphi}(x_2)\not\in \Pi.$
Set
$$t_1:=\inf\limits_{t\in [0, 1], \widetilde{\varphi}
(\widetilde{\gamma}(t))\in \Pi}t\,,\qquad t_2:=
\sup\limits_{t\in [0, 1],
\widetilde{\varphi}(\widetilde{\gamma}(t))\in \Pi}t\,.$$
Since by the assumption $|\widetilde{\gamma}|\cap |\alpha_1|\ne
\varnothing,$ by~\cite[Theorem~1.I.5.46]{Ku} we obtain that
$\widetilde{\varphi}(\widetilde{\gamma}(t_1))\in \partial \Pi$ and
$\widetilde{\varphi}(\widetilde{\gamma}(t_2))\in \partial \Pi.$ Now,
we may replace a path $\widetilde{\gamma}|_{[t_1, t_2]}$ by a path
$\alpha_*:[t_1, t_2]\rightarrow {\Bbb R}^2$ which does not intersect
$I.$ Finally, set
$$\gamma(t)=\begin{cases}\widetilde{\gamma}(t)\,, &t\in[0, 1]\setminus[t_1, t_2]\,, \\
\widetilde{\varphi}^{\,-1}(\alpha_*(t))\,,& t\in [t_1,
t_2]\end{cases}\,.$$
The path $\gamma$ joins any $x_1, x_2\in {\Bbb R}^n\setminus
(\overline{B(z_0, \varepsilon_2)}\cup |\alpha_1|)$ by a path
$\gamma:[0, 1]\rightarrow {\Bbb R}^n,$ $\gamma(0)=x_1,$
$\gamma(1)=x_2,$ in ${\Bbb R}^n\setminus |\alpha_1|.$

\medskip
In any of two cases, $n=2$ or $n\geqslant 3,$ we have proved that,
we may join any $x_1, x_2\in {\Bbb R}^n\setminus (\overline{B(z_0,
\varepsilon_2)}\cup |\alpha_1|)$ by a path $\gamma:[0, 1]\rightarrow
{\Bbb R}^n,$ $\gamma(0)=x_1,$ $\gamma(1)=x_2,$ in ${\Bbb
R}^n\setminus |\alpha_1|.$ Let us show that $\gamma$ may be chosen
in ${\Bbb R}^n\setminus (\overline{B(z_0, \varepsilon_2)}\cup
|\alpha_1|),$ as well. Choose $\varepsilon_3>\varepsilon_2>0$ such
that $\overline{B(z_0, \varepsilon_3)}\cap |\alpha_1|=\varnothing.$
If $|\gamma|\cap \overline{B(z_0, \varepsilon_2)}\ne\varnothing,$ by
\cite[Theorem~1.I.5.46]{Ku} we have that $|\gamma|\cap S(z_0,
\varepsilon_3)\ne\varnothing.$ Let
$$t_1:=\inf\limits_{t\in [0, 1]}
\gamma(t)\in S(z_0, \varepsilon_3)\,,\qquad t_2:=\sup\limits_{t\in
[0, 1]} \gamma(t)\in S(z_0, \varepsilon_3)\,.$$
Since $S(z_0, \varepsilon_3)$ is connected, we may join the points
$\gamma(t_1)$ and $\gamma(t_2)$ in $S(z_0, \varepsilon_3)$ by some a
path $\alpha_{**}:[t_1, t_2]\rightarrow S(z_0, \varepsilon_3).$
Finally,
$$\widetilde{\gamma}(t)=\begin{cases}\gamma(t)\,, &t\in[0, 1]\setminus[t_1, t_2]\,, \\
\alpha_{**}(t)\,,& t\in [t_1, t_2]\end{cases}$$
is a required path, because $\widetilde{\gamma}$ joins $x_1$ and
$x_2$ in ${\Bbb R}^n\setminus (\overline{B(z_0, \varepsilon_2)}\cup
|\alpha_1|).$ Thus, ${\Bbb R}^n\setminus (\overline{B(z_0,
\varepsilon_2)}\cup |\alpha_1|)$ is a domain, i.e.,
$R=R(\overline{B(z_0, \varepsilon_2)}, |\alpha_1|)$ is a ring
domain.

By Proposition~\ref{pr2} and by~(\ref{eq2A})
$$M(\Gamma(\overline{B(z_0, \varepsilon_2)}, |\alpha_1|,
\overline{{\Bbb R}^n}))\geqslant\frac{\omega_{n-1}}{\left(
\log\frac{2\lambda^2_n}{h(\overline{B(z_0,
\varepsilon_2)})h(|\alpha_1|)}\right)^{n-1}}\geqslant$$
\begin{equation}\label{eq4A}
\geqslant \frac{\omega_{n-1}}{\left(
\log\frac{2\lambda^2_n}{h(\overline{B(z_0, \varepsilon_2)})c_1\cdot
|x-x_0|}\right)^{n-1}}\,,
\end{equation}
where $c_1:=\frac{1}{1+\varepsilon^2_1},$ moreover, $h(x,
y)\geqslant c_1|x-y|$ for any $x,y \in B(x_0, \varepsilon_1).$ Since
$D$ is a $QED$-domain, from~(\ref{eq4A}) it follows that there is
$c_2>0$ such that
\begin{equation}\label{eq5C}
M(\Gamma)\geqslant\frac{c_2\cdot\omega_{n-1}}{\left(
\log\frac{2\lambda^2_n}{h(\overline{B(z_0, \varepsilon_2)})c_1\cdot
|x-x_0|}\right)^{n-1}}\,,
\end{equation}
where $\Gamma:=\Gamma(\overline{B(z_0, \varepsilon_2)}, |\alpha_1|,
D).$ Observe that, (\ref{eq5C}) holds even if $\overline{B(z_0,
\varepsilon_2)}\cap |\alpha_1|\ne\varnothing,$ because the left part
of it equals to~$\infty.$

\medskip
Set%
$$\varepsilon:=|f(x)-f(x_0)|+|f(x_0)-y_0|\,.$$
On the other hand, we observe that
\begin{equation}\label{eq3G}
f(\Gamma)>\Gamma(S(y_0, \varepsilon), S(y_0, \varepsilon_0), A(y_0,
\varepsilon, \varepsilon_0))\,.
\end{equation}
Indeed, let $\widetilde{\gamma}\in f(\Gamma).$ Then
$\widetilde{\gamma}(t)=f(\gamma(t)),$ where $\gamma\in \Gamma,$
$\gamma:[0, 1]\rightarrow D,$ $\gamma(0)\in \overline{B(z_0,
\varepsilon_2)},$ $\gamma(1)\in |\alpha_{1}|.$ By the
relation~(\ref{eq1E}), we obtain that $f(\gamma(0))\in {\Bbb
R}^n\setminus B(y_0, \varepsilon_0),$ however, by the triangle
inequality and due to~(\ref{eq1E})
$$|w-y_0|\leqslant |w-f(x_0)|+|f(x_0)-f(x)|\leqslant
|f(x)-f(x_0)|+|f(x_0)-f(x)|<\varepsilon_0$$
for any $w\in |I(t)|,$ i.e., $|I|\subset B(y_0, \varepsilon_0).$
Thus,
$$f(\gamma(1))\subset f(|\alpha_1|)\subset |I|\subset B(y_0, \varepsilon_0)\,.$$
Therefore, $|f(\gamma(t))|\cap B(y_0, \varepsilon_0)\ne\varnothing
\ne |f(\gamma(t))|\cap ({\Bbb R}^n\setminus B(y_0, \varepsilon_0)).$
Now, by~\cite[Theorem~1.I.5.46]{Ku} we obtain that, there is
$0<t_1<1$ such that $f(\gamma(t_1))\in S(y_0, \varepsilon_0).$ Set
$\gamma_1:=\gamma|_{[t_1, 1]}.$ We may consider that
$f(\gamma(t))\in B(y_0, \varepsilon_0)$ for any $t\geqslant t_1.$
Further, $f(\gamma(0))\in {\Bbb R}^n\setminus B(y_0,
\varepsilon_0)\in {\Bbb R}^n\setminus B(y_0,\varepsilon),$ because
$B(y_0, \varepsilon)\subset B(y_0,\varepsilon_0)$ by the first
relation in~(\ref{eq1E}). On the other hand, by the triangle
inequality
$$|w-y_0|\leqslant |w-f(x_0)|+|f(x_0)-f(x)|=\varepsilon$$
for any $w\in |I(t)|.$ Thus, $|f(\gamma(t))|\cap \overline{B(y_0,
\varepsilon)}\ne\varnothing \ne |f(\gamma(t))|\cap ({\Bbb
R}^n\setminus \overline{B(y_0, \varepsilon)}).$
By~\cite[Theorem~1.I.5.46]{Ku} we obtain that, there is $t_2\in
[t_1, 1]$ such that $f(\gamma(t_2))\in S(y_0, \varepsilon).$ Put
$\gamma_2:=\gamma|_{[t_1, t_2]}.$ We may consider that
$f(\gamma(t))\not\in B(y_0, \varepsilon)$ for any $t\in [t_1, t_2].$
Now, the path $f(\gamma_2)$ is a subpath of
$f(\gamma)=\widetilde{\gamma},$ which belongs to $\Gamma(S(y_0,
\varepsilon), S(y_0, \varepsilon_0), A(y_0, \varepsilon,
\varepsilon_0)).$ The relation~(\ref{eq3G}) is established.

\medskip
It follows from~(\ref{eq3G}) that
\begin{equation}\label{eq3H}
\Gamma>\Gamma_{f}(S(y_0, \varepsilon), S(y_0, \varepsilon_0), A(y_0,
\varepsilon, \varepsilon_0))\,.
\end{equation}
By the assumption, $I(\varepsilon, \varepsilon_0)>0$ for all
$\varepsilon\in (0, \varepsilon_0).$ Set
$$\eta(t)=\left\{
\begin{array}{rr}
\psi(t)/I(|f(x)-f(x_0)|, \varepsilon_0), & t\in (\varepsilon, \varepsilon_0)\,,\\
0,  &  t\not\in (\varepsilon, \varepsilon_0)\,,
\end{array}
\right. $$
where $I(\varepsilon,
\varepsilon_0)=\int\limits_{\varepsilon}^{\varepsilon_0}\,\psi (t)\,
dt.$ Observe that
$\int\limits_{\varepsilon}^{\varepsilon_0}\eta(t)\,dt=1.$ Now, by
the relations~(\ref{eq3.7.2}) and~(\ref{eq3H}), and due to the
definition of $f$ in~(\ref{eq2*A})--(\ref{eqA2}), we obtain that
$$M(\Gamma)\leqslant M(\Gamma_{f}(S(y_0, \varepsilon), S(y_0,
\varepsilon_0), A(y_0, \varepsilon, \varepsilon_0)))\leqslant$$
\begin{equation}\label{eq3J}
\leqslant \frac{1}{I^n(\varepsilon,
\varepsilon_0)}\int\limits_{A(y_0, \varepsilon, \varepsilon_0)}
Q(y)\cdot\psi^{\,n}(|y-y_0|)\,dm(y)=\alpha(\varepsilon,
\varepsilon_0)\,.
\end{equation}
Combining~(\ref{eq5C}) with~(\ref{eq3J}), we obtain that
$$\frac{\omega_{n-1}c_2}{\left(
\log\frac{2\lambda^2_n}{h(\overline{B(z_0, \varepsilon_2)})c_1\cdot
|x-x_0|}\right)^{n-1}}\leqslant\alpha(\varepsilon,
\varepsilon_0)\,.$$
Expressing $|x-x_0|$ in this relation, we obtain the desired
relation~(\ref{eq1}). Lemma is proved.~$\Box$

\medskip
\begin{lemma}\label{lem4}
{\it Let $D$ be a domain in ${\Bbb R}^n,$ $n\geqslant 2,$ and let
$f_j:D\rightarrow {\Bbb R}^n,$ $n\geqslant 2,$ $j=1,2,\ldots,$ be a
homeomorphisms satisfying the conditions~(\ref{eq2*A})--(\ref{eqA2})
at a point $y_0\in {\Bbb R}^n$ and converging to some mapping
$f:D\rightarrow \overline{{\Bbb R}^n}$ as $j\rightarrow\infty$
locally uniformly in $D$ with respect to the chordal metric $h.$
Assume that, $f$ is not a constant in $D.$ Then for any $y_0\in
{\Bbb R}^n$ there is $\varepsilon_0=\varepsilon_0(y_0)>0,$ $z_0\in
D$ and $\varepsilon_2=\varepsilon_2(z_0)>0$ such that
\begin{equation}\label{eq4} f_m(E)\cap \overline{B(y_0,
\varepsilon_0)}=\varnothing,\qquad m=1,2,\ldots \,,
\end{equation}
where $E:=\overline{B(z_0, \varepsilon_2)}.$ }
\end{lemma}

\medskip
\begin{proof}
Since $f$ is not a constant in $D,$ there are $u, w\in B(x_0,
\varepsilon_1)$ such that $f(u)\ne f(v).$ By the convergence of
$f_m$ to $f,$ we have that
\begin{equation}\label{eq3A}
h(f_m(u), f_m(v))\geqslant \delta>0
\end{equation}
for some $\delta>0$ and all $m=1,2,\ldots. $

\medskip
Let $E_1$ be a path joining $u$ and $v$ in $D.$ Put
$0<\varepsilon_0=\varepsilon_0(y_0)<\delta/2.$ Since by~(\ref{eq3A})
$h(f_m(E_1))\geqslant \delta$ for any $m\in{\Bbb N}$ and
$d(f_m(E_1))\geqslant h(f_m(E_1)),$
\begin{equation}\label{eq2}
f_m(E_1)\setminus B(y_0, \varepsilon_0)\ne \varnothing,\qquad
m=1,2,\ldots \,.
\end{equation}
By~(\ref{eq2}), there is $w_m=f_m(z_m)\in \overline{{\Bbb
R}^n}\setminus B(y_0, \varepsilon_0),$ where $z_m\in E_1.$ Since
$E_1$ is a continuum, $\overline{{\Bbb R}^n}$ is a compactum and the
set $\overline{{\Bbb R}^n}\setminus B(y_0, \varepsilon_0)$ is
closed, we may consider that $z_m\rightarrow z_0\in E_1$ as
$m\rightarrow\infty$ and $w_m\rightarrow w_0\in \overline{{\Bbb
R}^n}\setminus B(y_0, \varepsilon_0).$ Obviously, $w_0\ne y_0.$

\medskip
Since $f_m$ converges to $f$ locally uniformly, the family $f_m$ is
equicontinuous due to Arzela-Ascoli theorem (see, e.g.,
\cite[item~20.4]{Va}). Thus, for any $\sigma>0$ there is
$\varepsilon_2=\varepsilon_2(z_0)>0$ such that $h(f_m(z_0),
f_m(z))<\sigma$ whenever $|z-z_0|\leqslant \varepsilon_2.$ Then, by
the triangle inequality
\begin{equation}\label{eq3}
h(f_m(z), w_0)\leqslant h(f_m(z), f_m(z_0))+ h(f_m(z_0), f_m(z_m))+
h(f_m(z_m), w_0)<3\sigma
\end{equation}
for $|z-z_0|\leqslant\varepsilon_2,$ some $M_1\in {\Bbb N}$ and all
$m\geqslant M_1.$ We may consider that latter holds for any
$m=1,2,\ldots .$ Since $w_0\in \overline{{\Bbb R}^n}\setminus B(y_0,
\varepsilon_0),$ we may choose $\sigma>0$ such that
$\overline{B_h(w_0, 3\sigma)}\cap \overline{B(y_0,
\varepsilon_2)}=\varnothing,$ where $B_h(w_0, \sigma)=\{w\in
\overline{{\Bbb R}^n}: h(w, w_0)<\sigma\}.$ Then~(\ref{eq3}) implies
that
\begin{equation}\label{eq4AA}
f_m(E)\cap \overline{B(y_0, \varepsilon_0)}=\varnothing,\qquad
m=1,2,\ldots \,,
\end{equation}
where $E:=\overline{B(z_0, \varepsilon_2)},$ as required.~$\Box$
\end{proof}

\medskip
\begin{lemma}\label{lem3}
{\it Let $D$ be a domain in ${\Bbb R}^n,$ $n\geqslant 2,$ and let
$f_j:D\rightarrow {\Bbb R}^n,$ $n\geqslant 2,$ $j=1,2,\ldots,$ be a
homeomorphisms satisfying the conditions~(\ref{eq2*A})--(\ref{eqA2})
at any point $y_0\in \overline{{\Bbb R}^n}$ and converging to some
mapping $f:D\rightarrow \overline{{\Bbb R}^n}$ as
$j\rightarrow\infty$ locally uniformly in $D$ with respect to the
chordal metric $h.$ Let $x_0\in D$ and let $B(x_0,
\varepsilon_1)\subset D$ such that $f$ is not a constant in $B(x_0,
\varepsilon_1).$ Assume that, for any $y_0\in \overline{{\Bbb R}^n}$
there is $\varepsilon_0=\varepsilon_0(y_0)>0$ and a Lebesgue
measurable function $\psi:(0, \varepsilon_0)\rightarrow [0,\infty]$
such that
\begin{equation}\label{eq9} I(\varepsilon,
\varepsilon_0):=\int\limits_{\varepsilon}^{\varepsilon_0}\psi(t)\,dt
< \infty\quad \forall\,\,\varepsilon\in (0, \varepsilon_0)\,,\quad
I(\varepsilon, \varepsilon_0)\rightarrow
\infty\quad\text{as}\quad\varepsilon\rightarrow 0\,,
\end{equation}
and, in addition, for some $\alpha=\alpha(\varepsilon,
\varepsilon_0)>0$
\begin{equation} \label{eq10}
\int\limits_{A(y_0, \varepsilon, \varepsilon_0)}
Q(y)\cdot\psi^{\,n}(|y-y_0|)\,dm(x)=\alpha(\varepsilon,
\varepsilon_0)\cdot I^n(\varepsilon, \varepsilon_0)\,,\end{equation}
as $\varepsilon\rightarrow 0,$ where $A(y_0, \varepsilon,
\varepsilon_0)$ is defined in (\ref{eq1**}). Assume that
$\alpha(\varepsilon, \varepsilon_0)\rightarrow 0$ as
$\varepsilon\rightarrow 0.$ Then there is $r_0>0$ such that
\begin{equation}\label{eq1A}
f_m(B(x_0, \varepsilon_1))\supset B(f_m(x_0), r_0)\,,\qquad
m=1,2,\ldots\,.
\end{equation}
}
\end{lemma}
\begin{remark}
If $y_0=\infty,$ the relation~(\ref{eq3.7.2}) must be understood by
the using the inversion $\psi_1(y)=\frac{y}{|y|^2}$ at the origin.
In other words, instead of
$$\int\limits_{A(y_0, \varepsilon, \varepsilon_0)}
Q(y)\cdot\psi^{\,n}(|y-y_0|)\,dm(y) = \alpha(\varepsilon,
\varepsilon_0)\cdot I^n(\varepsilon, \varepsilon_0)$$
we need to consider the condition
$$\int\limits_{A(0, \varepsilon, \varepsilon_0)}
Q\left(\frac{y}{|y|^2}\right)\cdot\psi^{\,n}(|y|)\,dm(y) =
\alpha(\varepsilon, \varepsilon_0)\cdot I^n(\varepsilon,
\varepsilon_0)\,.$$
\end{remark}

\medskip {\it Proof of Lemma~\ref{lem3}.} Assume the contrary. Then
there is $r_m>0,$ $m=1,2,\ldots\,,$ $r_m\rightarrow 0$ as
$m\rightarrow\infty,$ and $y_m\in B(f_m(x_0), r_m)$ such that
$y_m\not\in f_m(B(x_0, \varepsilon_1)).$ Join $y_m$ and $f_m(x_0)$
by a path $\gamma_m$ inside $B(f_m(x_0), r_m).$ Let $\alpha_m:[0,
c)\rightarrow B(x_0, \varepsilon_1)$ be a maximal $f_m$-lifting of
$\gamma_m$ starting at $x_0$ in $B(x_0, \varepsilon_1).$ This
lifting exists by Proposition~\ref{pr3}; by the same Proposition we
have that $\alpha_m\rightarrow S(x_0, \varepsilon_1)$ as
$t\rightarrow c-0.$ So, we may find $\omega_m\in |\alpha_m|\subset
B(x_0, \varepsilon_1)$ such that $d(\omega_m, S(x_0,
\varepsilon_1))<1/m.$ We may assume that $\omega_m\rightarrow
\omega_0\in S(x_0, \varepsilon_1).$

\medskip Since $f_m$ converges to $f$ locally uniformly,
$f_m(x_0)\rightarrow f(x_0)$ as $m\rightarrow \infty.$ Without loss
of generality, we may assume that $f(x_0)\ne\infty;$ in other case
we consider the family $\psi_1\circ f_m$ instead of $f_m,$ where
$\psi_1(x)=\frac{x}{|x|^2}.$ Set $y_0=f(x_0).$

\medskip
By Lemma~\ref{lem4}, there is $\widetilde{\varepsilon_0}>0,$ $z_0\in
D$ and $\varepsilon_2=\varepsilon_2(z_0)>0$ such that $f_m(E)\cap
\overline{B(y_0, \widetilde{\varepsilon_0})}=\varnothing,$ $
m=1,2,\ldots \,,$ where $E:=\overline{B(z_0, \varepsilon_2)}.$ We
may consider that $\varepsilon_0<\widetilde{\varepsilon_0},$ where
$\varepsilon_0$ is from~(\ref{eq9})--(\ref{eq10}). Thus,
\begin{equation}\label{eq4B} f_m(E)\cap \overline{B(y_0,
\varepsilon_0})=\varnothing,\qquad m=1,2,\ldots \,.
\end{equation}
Since $I(\varepsilon, \varepsilon_0)\rightarrow \infty$ as
$\varepsilon\rightarrow 0$ and $I(\varepsilon,
\varepsilon_0)<\infty,$ we have that $I(\varepsilon,
\varepsilon_0)>0$ for sufficiently small $\varepsilon.$ Observe
that, $B(x_0, \varepsilon_1)$ is a $QED$-domain,
see~\cite[Lemma~4.3]{Vu}.

We may apply Lemma~\ref{lem1} for $x:=\omega_m$ and
$\varepsilon=\varepsilon_m:=|f_m(\omega_m)-f_m(x_0)|+|f_m(x_0)-y_0|.$
Observe that, $|f_m(\omega_m)-f_m(x_0)|<r_m,$ $r_m\rightarrow 0,$
because by the construction $f_m(\omega_m)\subset |\alpha_m|\subset
B(f_m(x_0), r_m).$ Thus, $\varepsilon_m\rightarrow 0$ as
$m\rightarrow \infty.$ Since $I(\varepsilon,
\varepsilon_0)\rightarrow \infty$ as $\varepsilon\rightarrow 0,$ the
relation~(\ref{eq1}) together with~(\ref{eq9})--(\ref{eq10}) yields
\begin{equation}\label{eq2B}
|\omega_m-x_0|\leqslant \frac{2\lambda^2_n}{c_1\cdot
h(\overline{B(z_0,
\varepsilon_2)})}\cdot\exp\left\{-\frac{\omega_{n-1}c_2}{\alpha(\varepsilon_m,
\varepsilon_0)}\right\}\rightarrow 0\,,\qquad m\rightarrow\infty\,,
\end{equation}
which is impossible because by the construction $\omega_m\rightarrow
\omega_0\in S(x_0, \varepsilon_1)$ as $m\rightarrow\infty,$ so
$|\omega_m-x_0|\geqslant \delta_*>0$ for sufficiently large
$m=1,2,\ldots .$ The contradiction obtained above
proves~(\ref{eq1A}).~$\Box$

\medskip
\begin{lemma}\label{lem2}
{\it Let $D$ be domain in ${\Bbb R}^n,$ $n\geqslant 2,$ and let
$f_j:D\rightarrow {\Bbb R}^n,$ $n\geqslant 2,$ $j=1,2,\ldots,$ be a
homeomorphisms satisfying the conditions~(\ref{eq2*A})--(\ref{eqA2})
at any point $y_0\in \overline{{\Bbb R}^n}$ and converging to some
mapping $f:D\rightarrow \overline{{\Bbb R}^n}$ as
$j\rightarrow\infty$ locally uniformly in $D$ with respect to the
chordal metric $h.$ Assume that, for any $y_0\in \overline{{\Bbb
R}^n}$ there is $\varepsilon_0=\varepsilon_0(y_0)>0$ and a Lebesgue
measurable function $\psi:(0, \varepsilon_0)\rightarrow [0,\infty]$
such that
\begin{equation}\label{eq7} I(\varepsilon,
\varepsilon_0):=\int\limits_{\varepsilon}^{\varepsilon_0}\psi(t)\,dt
< \infty\quad \forall\,\,\varepsilon\in (0, \varepsilon_0)\,,\quad
I(\varepsilon, \varepsilon_0)\rightarrow
\infty\quad\text{as}\quad\varepsilon\rightarrow 0\,,
\end{equation}
and, in addition, for some $\alpha(\varepsilon, \varepsilon_0)>0$
\begin{equation} \label{eq8}
\int\limits_{A(y_0, \varepsilon, \varepsilon_0)}
Q(y)\cdot\psi^{\,n}(|y-y_0|)\,dm(x)=\alpha(\varepsilon,
\varepsilon_0)\cdot I^n(\varepsilon, \varepsilon_0)\,,\end{equation}
as $\varepsilon\rightarrow 0,$ where $A(y_0, \varepsilon,
\varepsilon_0)$ is defined in (\ref{eq1**}). Assume that
$\alpha(\varepsilon, \varepsilon_0)\rightarrow 0$ as
$\varepsilon\rightarrow 0.$ Then $f$ is discrete.}
\end{lemma}

\medskip
\begin{proof}
Assume the contrary. Then there is $x_0 \in D$ and a sequence
$x_m\in D,$ $m=1,2,\ldots ,$ $x_m \ne x_0,$ such that
$x_m\rightarrow x_0$ as $m\rightarrow \infty$ and $f(x_m)=f(x_0).$
Observe that, $E_0=\{x\in D: f(x)=f(x_0)\}$ is closed in $D$ by the
continuity of $f$ and does not coincide with $D,$ because $f\not
\equiv const.$ Thus, we may consider that $x_0$ may be replaced by
non isolated boundary point of $E_0.$

\medskip
Let us prove that $f(x_0)\ne \infty$ for any $x_0\in D.$ Let $x_0\in
D$ and let $y_0=f(x_0).$ Since $f$ is not a constant in $B(x_0,
\varepsilon_1),$ by Lemma~\ref{lem3} there is $r_0>0,$ which does
not depend on $m,$ such that
$B(f_m(x_0), r_0)\subset f_m(B(x_0, \varepsilon_1)),$
$m=1,2,\ldots.$ Then also $B_h(f_m(x_0), r_*)\subset f_m(B(x_0,
\varepsilon_1)),$ $m=1,2,\ldots,$ for some $r_*>0.$ Let $y\in
B_h(y_0, r_*/2)=B_h(f(x_0), r_*/2).$ By the converges of $f_m$ to
$f$ and by the triangle inequality, we obtain that
$$h(y, f_m(x_0))\leqslant h(y, f(x_0))+h(f(x_0), f_m(x_0))<r_*/2+r_*/2=r_*$$
for sufficiently large $m\in {\Bbb N}.$ Thus,
$$B_h(f(x_0), r_*/2)\subset B_h(f_m(x_0), r_*)\subset f_m(B(x_0,
\varepsilon_1))\subset {\Bbb R}^n\,.$$
In particular, $y_0=f(x_0)\in {\Bbb R}^n,$ as required.

\medskip
Now, $f:U\rightarrow {\Bbb R}^n,$ where $U$ is a some neighborhood
of $x_0.$ Since $f_m$ converges to $f$ locally uniformly,
$f_m(x_0)\rightarrow f(x_0)$ as $m\rightarrow \infty.$ By the
proving above, $f(x_0)\ne\infty;$ thus, $f_m$ converges uniformly to
$f$ in $U$ by the Euclidean metric, as well. We may consider that
$U$ is a domain.

By Lemma~\ref{lem4} we may construct a continuum $E\subset D$ for
which~(\ref{eq4AA}) holds, where $\varepsilon_0=\varepsilon_0(y_0)$
is some number. Decreasing $\varepsilon_0,$ if it is required, we
may consider that $\varepsilon_0$ is a number
from~(\ref{eq7})--(\ref{eq8}). Observe that, the inequality
$$\varepsilon=\varepsilon_{j,
m}:=|f_j(x_m)-f_j(x_0)|+|f_j(x_0)-y_0|<\varepsilon_0$$
holds for all $j\geqslant j_0=j_0(\varepsilon_0)$ and $m\geqslant
m_0=m_0(\varepsilon_2)$ due to the local uniform convergence of
$f_j$ (and thus by the equicontinuity of the family $f_m,$
$m=1,2,\ldots $). Applying Lemma~\ref{lem1}, we obtain that
$$|x_m-x_0|\leqslant $$
\begin{equation}\label{eq2C}
\leqslant \frac{2\lambda^2_n}{c_1\cdot h(\overline{B(z_0,
\varepsilon_2)})}\cdot\exp\left\{-\frac{\omega_{n-1}c_2}
{\alpha(|f_j(x_m)-f_j(x_0)|+|f_j(x_0)-y_0|,
\varepsilon_0)}\right\}\,.
\end{equation}
Taking here the limit as $j\rightarrow \infty,$ we obtain that
$$|x_m-x_0|\leqslant $$
\begin{equation}\label{eq2D}
\leqslant \frac{2\lambda^2_n}{c_1\cdot h(\overline{B(z_0,
\varepsilon_2)})}\cdot\exp\left\{-\frac{\omega_{n-1}c_2}{\alpha(|f(x_m)-f(x_0)|,
\varepsilon_0)}\right\}\,.
\end{equation}
By~(\ref{eq2D}) $f(x_m)\ne f(x_0)$ for $m\geqslant
m_0=m_0(\varepsilon_2).$ Indeed, in the contrary case the right hand
side of~(\ref{eq2D}) must equals to $0$ for $m\geqslant
m_0=m_0(\varepsilon_2)$ due to the condition $I(\varepsilon,
\varepsilon_0)\rightarrow\infty.$ However, this contradicts to the
inequality~(\ref{eq2D}) because $|x_m-x_0|>0$ by the choice of
$x_m.$

Finally, $f$ is discrete, because the relation $f(x_m)\ne f(x_0)$
for $m\geqslant m_0=m_0(\varepsilon_2)$ contradicts the assumption
made above.~$\Box$
\end{proof}

\medskip
\begin{lemma}\label{lem5}
{\it Under assumptions of Lemma~\ref{lem2}, either $f$ is a constant
in $\overline{{\Bbb R}^n}$, or $f$ is a homeomorphism
$f:D\rightarrow{\Bbb R}^n.$}
\end{lemma}

\medskip
\begin{proof}
Let $f$ is not a constant. By Lemma~\ref{lem2}, $f:D\rightarrow
\overline{{\Bbb R}^n}$ is discrete. By Proposition~\ref{pr4}
$f:D\rightarrow \overline{{\Bbb R}^n}$ is a homeomorphism. Let
$x_0\in D$ and $y_0=f(x_0).$ We put $\varepsilon_1>0$ such that
$B(x_0, \varepsilon_1)\subset D.$ By Lemma~\ref{lem3} there is
$r_0>0$ such that
\begin{equation}\label{eq1B}
f_m(B(x_0, \varepsilon_1))\supset B(f_m(x_0), r_0)\,,\qquad
m=1,2,\ldots\,.
\end{equation}
Then also $B_h(f_m(x_0), r_*)\subset f_m(B(x_0, \varepsilon_1)),$
$m=1,2,\ldots,$ for some $r_*>0.$ Let $y\in B_h(y_0,
r_*/2)=B_h(f(x_0), r_*/2).$ By the converges of $f_m$ to $f$ and by
the triangle inequality, we obtain that
$$h(y, f_m(x_0))\leqslant h(y, f(x_0))+h(f(x_0), f_m(x_0))<r_*/2+r_*/2=r_*$$
for sufficiently large $m\in {\Bbb N}.$ Thus,
$$B_h(f(x_0), r_*/2)\subset B_h(f_m(x_0), r_*)\subset f_m(B(x_0,
\varepsilon_1))\subset {\Bbb R}^n\,.$$
In particular, $y_0=f(x_0)\in {\Bbb R}^n,$ as required.
\end{proof}

\section{Proof of the main result}

The following statement may be found in~\cite[Lemma~1.3]{Sev}.

\medskip
\begin{proposition}\label{pr6}
{\it\, Let $Q:{\Bbb R}^n\rightarrow [0,\infty],$ $n\geqslant 2,$ be
a Lebesgue measurable function and let $x_0\in {\Bbb R}^n.$ Assume
that either of the following conditions holds

\noindent (a) $Q\in FMO(x_0),$

\noindent (b)
$q_{x_0}(r)\,=\,O\left(\left[\log{\frac1r}\right]^{n-1}\right)$ as
$r\rightarrow 0,$

\noindent (c) for some small $\delta_0=\delta_0(x_0)>0$ we have the
relations
\begin{equation}\label{eq5***B}
\int\limits_{\delta}^{\delta_0}\frac{dt}{tq_{x_0}^{\frac{1}{n-1}}(t)}<\infty,\qquad
0<\delta<\delta_0,
\end{equation}
and
\begin{equation}\label{eq5**}
\int\limits_{0}^{\delta_0}\frac{dt}{tq_{x_0}^{\frac{1}{n-1}}(t)}=\infty\,.
\end{equation}
Then  there exist a number $\varepsilon_0\in(0,1)$ and a function
$\psi(t)\geqslant 0$ such that the relation
\begin{equation}\label{eqlem21}
\int\limits_{\varepsilon<|x-b|<\varepsilon_0}Q(x)\cdot\psi^n(|x-b|)
\ dm(x)=o(I^n(\varepsilon, \varepsilon_0))\,,
\end{equation}
holds as $\varepsilon\rightarrow 0,$
where $\psi:(0, \varepsilon_0)\rightarrow [0, \infty)$ is some
function such that, for some $0<\varepsilon_1<\varepsilon_0,$
\begin{equation}\label{eqlem22}
0<I(\varepsilon, \varepsilon_0)
=\int\limits_{\varepsilon}^{\varepsilon_0}\psi(t)\,dt < \infty
\qquad \forall\quad\varepsilon \in(0, \varepsilon_1)\,.
\end{equation}}
\end{proposition}

\medskip
{\it Proof of Theorem~\ref{th1}} immediately follows by
Lemma~\ref{lem5} and Proposition~\ref{pr6}.~$\Box$

%=================Список литературы====================
%\end{fulltext}

CONTACT INFORMATION

\medskip
{\bf \noindent Evgeny Sevost'yanov} \\
{\bf 1.} Zhytomyr Ivan Franko State University,  \\
40 Bol'shaya Berdichevskaya Str., 10 008  Zhytomyr, UKRAINE \\
{\bf 2.} Institute of Applied Mathematics and Mechanics\\
of NAS of Ukraine, \\
19 Henerala Batyuka Str., 84 116 Slavyansk,  UKRAINE\\
esevostyanov2009@gmail.com

\medskip
{\bf \noindent Valery Targonskii} \\
Zhytomyr Ivan Franko State University,  \\
40 Bol'shaya Berdichevskaya Str., 10 008  Zhytomyr, UKRAINE \\
w.targonsk@gmail.com

\end{document}